\documentclass[pdflatex,sn-mathphys-num]{sn-jnl}% Math and Physical Sciences Numbered Reference Style
\usepackage{graphicx}%
\usepackage{multirow}%
\usepackage{amsmath,amssymb,amsfonts}%
\usepackage{amsthm}%
\usepackage{mathrsfs}%
\usepackage[title]{appendix}%
\usepackage{xcolor}%
\usepackage{textcomp}%
\usepackage{manyfoot}%
\usepackage{booktabs}%
\usepackage{algorithm}%
\usepackage{algorithmicx}%
\usepackage{algpseudocode}%
\usepackage{listings}%
\theoremstyle{thmstyleone}%
\newtheorem{theorem}{Theorem}%  meant for continuous numbers
\theoremstyle{thmstyletwo}%
\newtheorem{remark}{Remark}%
\newtheorem{lemma}{Lemma}
\theoremstyle{thmstylethree}%

\def \Rn {\mathbb{R}^n}

\begin{document}

        \title[A Remark on certain Oscillatory Singular Integrals with Nonstandard   Kernel ]{A New Proof About Certain  Oscillatory Singular Integrals with Nonstandard   Kernel }

%%  \sfx{IV}}\email{iauthormail.com}
%%=============================================================%%

\author*[1]{\fnm{Shen} \sur{ Jiawei}}\email{sjiawei633@xtu.edu.cn}

\affil*[1]{\orgdiv{School of Mathematics and Computational Science}, \orgname{Xiangtan University}, \orgaddress{ \city{Xiangtan}, \postcode{411105}, \state{Hunan}, \country{China}}}


\abstract{In the paper, we provide  a new method to study the  oscillatory  singular  integral operator $T_{Q,A}$ with nonstandard kernel defined by
 \[T_{Q,A} f(x)=\text { p.v. } \int_{\mathbb{R}^{n}} f(y) \frac{\Omega(x-y)}{|x-y|^{n+1}}\left(A(x)-A(y)-\nabla A(y)(x-y)\right)  e^{i Q(|x-y|)} d y, \]
 where $Q(t)=\sum_{1\le i\le m} a_it^{\alpha_i}(a_i\in\mathbb{R} \text{and } a_i\neq 0, \alpha_i\in \mathbb{N})$ , and $\Omega$ is a homogeneous
function of degree zero on $\mathbb{R}^{n}$ and  satisfies the vanishing moment condition. Under the condition that  $\Omega\in L(logL)^2(\mathbb{S}^{n-1})$ and $\nabla A\in \text{BMO}(\mathbb{R}^n),$  the authors show that 
$T_{Q,A}$ is bounded on $L^p(\mathbb{R}^{n})$ with a uniform boundedness, which improves and extends the previous results. 
} 



\maketitle

\section{Introduction}
Let $n \geq 2$ and the unit sphere $S^{n-1}$ be equipped with the spherical Lebesgue measure d$\sigma$. For $x\in \Rn,$ we denote $x_j(1\le j\le n)$ as the $j$th variable of x and $x'=\frac{x}{|x|}(x\neq 0).$ In the paper, we assume that $\Omega \in L^{1}(\mathbb{S}^{n-1})$ is a homogeneous function of degree zero, and satisfies the vanishing moment condition that for all $1\le j\le n,$ 
\begin{equation}\label{2}
    \int_{S^{n-1}} \Omega\left(x^{\prime}\right) x'_j \text{d} \sigma\left(x^{\prime}\right)=0.
\end{equation}

Now give a suitable  function $A:\mathbb{R}^{n}\rightarrow \mathbb{R},$ we consider the following oscillatory  singular  integral operator $\tilde{T}_{P,A}$ defined by 
\[\tilde{T}_{P,A} f(x)=\text { p.v. } \int_{\mathbb{R}^{n}}  \frac{\Omega(x-y)}{|x-y|^{n+1}}(A(x)-A(y)-\nabla A(y)\cdot(x-y)) f(y)  e^{iP(x,y)} d y. \]
    About the operator $\tilde{T}_{P,A}$, when $P(x,y)$ is a constant,  we view it  as     \[T_{\Omega,A} f(x)=\text { p.v. } \int_{\mathbb{R}^{n}} \frac{\Omega(x-y)}{|x-y|^{n+1}}(A(x)-A(y)-\nabla A(y)\cdot(x-y)) f(y) d y.\]
  which is closely related to the $n-$dimensional Calder\'{o}n commutator, and,  when $\nabla A \in L^\infty(\Rn),$  $T_{\Omega,A}$ is viewed as the $n-$dimensional Calder\'{o}n commutator\cite{Calderon1965}.  However, when $\nabla A \in BMO(\Rn),$ $T_{\Omega,A}$ is  not a Calder\'{o}n-Zygmund operator even if  $\Omega\in \text{Lip}(S^{n-1}).$  When $\Omega\in \text{Lip}_\alpha(S^{n-1})(0<\alpha\le1),$ Cohen\cite{Cohen1981} pointed out that  $T_{\Omega,A}$ is bounded on $L^p(\Rn)$  for  $1<p<\infty.$  Hofmann \cite{Hofmann} improved the result of Cohen\cite{Cohen1981}, and showed that $\Omega\in L^\infty(\Rn)$ is sufficient condition such that $T_{\Omega,A}$ is bounded on $L^p(\Rn).$ Recently, many authors studied the boundness of  $T_{\Omega,A}$  when $\Omega$ satisfies certain minimum size condition.  Hu, Tao, Wang and Xue \cite{Hu2024} obtain the following result.

%   we can was first studied by Cohen\cite{Cohen1981}.    $\mathcal{C}_{\Omega,a}$, which is defined by 
%      \begin{equation}
% \mathcal{C}_{\Omega,a}f(x) = \text{p.v.} \int_{\mathbb{R}^d} \frac{\Omega(x - y)}{|x - y|^{n+ 1}} (a(x) - a(y)) f(y) \, dy
% \end{equation} 
% where $a$ is a function on $\Rn$ such that  $\partial_n a\in L^\infty(\Rn)$ for all $n$ with $1\le n\le d.$ This operator was introduced by Calder\'{o}n\cite{Calderon1965}.  After that,  many authors studied it and reader can refer to  \cite{calderon1,Calderon1965,calderon2,calderon3}, \cite[chapter 8]{Grafakos2}, \cite[chapter 4]{Muscalu}. Furthermore, Cohen\cite{Cohen1981} first studied the following operator $T_{\Omega,A}.$ 
                
% where $\nabla A\in BMO(\Rn).$    The operator $T_{\Omega,A}$ is closely related to the  $n-$dimensional Calder\'{o}n commutator $\mathcal{C}_{\Omega,a}$.   When $\nabla A \in L^\infty(\Rn),$  $T_{\Omega,A}$ is viewed as $n-$dimensional Calder\'{o}n commutator $\mathcal{C}_{\Omega,a}$. However,  $T_{\Omega,A}$ is  not a Calder\'{o}n-Zygmund operator even if  $\Omega\in \text{Lip}(S^{n-1}).$  When $\Omega\in \text{Lip}_\alpha(S^{n-1}),$ Cohen point out that  $T_{\Omega,A}$ is $L^p(\Rn)$ bounded for $1<p<\infty.$  Hofman improve the result of Cohen, and show that $\Omega\in L^\infty$ is sufficient condition such that $T_{\Omega,A}$ is bounded on $L^p(\Rn).$ Recently, many authors study the boundness of  $T_{\Omega,A}$  when $\Omega$ satisfies certain minimum size condition.  Hu, Tao, Wang and Xue \cite{Hu2024} pointed out that  

\begin{theorem}\label{TBOUND}
    Let $\Omega$ be homogeneous of degree zero and  satisfy the vanishing condition (1.1),  and $A$ be a function on $\Rn$ such that $\nabla A\in \text{BMO}(\mathbb{R}^n)$. Suppose that $\Omega\in L(\log L)^2(S^{n - 1})$. Then, $T_{\Omega,A}$ is bounded on $L^p(\mathbb{R}^n)$ for all $p\in(1,\infty).$ 
\end{theorem}
At the same time, Tao and Hu \cite{Tao2022} considered the maximal   operator $T_{\Omega, A}^{*}$ defined by  
           \[T_{\Omega, A}^{*} f(x): =\sup _{\epsilon>0}\left|\int_{|x-y|>\epsilon} \frac{\Omega(x-y)}{|x-y|^{n+1}}(A(x)-A(y)-\nabla A(y)(x-y)) f(y) d y\right|, \]
 and they got the following conclusion.

\begin{theorem}\label{MTbound}
     Let $\Omega$ be homogeneous of degree zero and satisfy the vanishing moment (1.1), and let $A$ be a function on $\mathbb{R}^n$ with derivatives of order one in $\text{BMO}(\mathbb{R}^n)$. Suppose that $\Omega\in L(\log L)^2(S^{n - 1}).$ Then $T^*_{\Omega,A}$ is  bounded on $L^p(\mathbb{R}^n)$ for all $p\in(1,\infty).$ 
\end{theorem}
On the other hand, when $P(x,y)$ is a polynomial about $x,y,$ many authors studied  it, and  readers can find  some results in \cite{ChenLu,Ding1996,Tao2022,Wang2025}.  Furthermore, in  \cite{ChenLu}, 
it is valuable to point out that the operator  $\tilde{T}_{P,A}$ is  ascribed to $T_{\Omega,A}$ when $\Omega\in L(\log L)^2(S^{n - 1}).$   So, it is not necessary to discuss the operator $ \tilde{T}_{P,A}$ .

In this paper,  we consider that  $Q(t)$ is a polynomial on $\mathbb{R}$ and we define the operator  $T_{Q,A}$ by
\begin{equation}\label{TQA}
    T_{Q,A} f(x)=\text { p.v. } \int_{\mathbb{R}^{n+1}} f(y) \frac{\Omega(x-y)}{|x-y|^{n+1}}\left(A(x)-A(y)-\nabla A(y)(x-y)\right)  e^{i Q(|x-y|)} d y,
    \end{equation}
which is  different from the case that  $P(x,y)$ is a polynomial on $\mathbb{R}^{2n}.$  Now, we get the following results:
% \begin{array}{l}
% L^{r}\left(\mathrm{~S}^{n-1}\right) \subsetneq L\left(\log L\right)^{\beta_{1}}\left(\mathrm{~S}^{n-1}\right) \\
% \subsetneq L\left(\log  L\right)^{\beta_{2}}\left(\mathrm{~S}^{n-1}\right) \\
% \forall r>1,0<\beta_{2}<\beta_{1},
% \end{array}
% $$
% and 

%  \[H^{1}\left(\mathrm{~S}^{n-1}\right) \subset  L\left(\log  L\right)^{\beta}\left(\mathrm{S}^{n-1}\right) \quad \forall \beta \geq 1 .\]
\begin{theorem}\label{main_th}
   Let $\Omega$ be homogeneous of degree zero and satisfy the vanishing moment \eqref{2}, and let $A$ be a function on $\mathbb{R}^n$ with derivatives of order one in $\text{BMO}(\mathbb{R}^n)$.   Let $m \in \mathbb{N}$, $Q(t)=\sum_{1\le i\le m} a_it^{\alpha_i}(a_i\in\mathbb{R}  \text{ and } a_i\neq 0, \alpha_i\in \mathbb{N})$, and  $T_{Q,A}$ be given in \eqref{TQA}. Suppose that $\Omega\in L(\log L)^2(\mathbb{S}^{n - 1})$, then, for all $1<p<\infty$,  $T_{Q,A}$ is bounded on $L^p$.That is,
    \begin{equation}\label{3}
        \left\|T_{Q,A} f\right\|_{L^{p}(\mathbb{R}^n)} \leq C_{m,A}(\|\Omega\|_{L(log L)^{2}(\mathbb{S}^{n - 1})}+1)\|f\|_{L^{p}(\mathbb{R}^n)}, 
    \end{equation}   
where $C_{m,A}$ depends only on   $m$, $A.$   
%     \begin{equation}\label{3}
%         \left\|T_{Q,A} f\right\|_{L^{p}(\mathbb{R}^n)(\mathbb{R}^n)} \leq C_m\|\Omega\|_{L(log L)^{2}(S^{d-1})}\|f\|_{L^{p}(\mathbb{R}^n)}, 
%     \end{equation}   
% where $C_m$ depends only on m, the number of monomials in $Q,$   not on the degree of $Q.$
\end{theorem}

 Our idea originates from Guo's work \cite{Guo2017}.  Firstly, we simply review Guo's proof method, and point out the difficulties  we need to overcome. In Guo's work, he firstly gave a fixed $\lambda$ and divided $\mathbb{R}$ into good scale and bad scale by Van der Corput lemma(Lemma \ref{15}). For bad scale, it is easy to use the Hardy-Litterwood Maximal function to console it. For good  scale, he found a lager scale $\lambda_k$ and $\gamma_k$   to divide the good scale into two parts : lower part and higher part such that $\ell'\le \gamma_k$ and $\ell'\ge \gamma_k.$ For lower part, it is easy to check the $L^p-$boundness.    In  higher  part,  $\ell'\ge \gamma_k,$ he use Plancherel’s Lemma and the Fourier decay estimate of the kernel  to get the follow inequality 
  \[
  \|H_{l'}^{(j_1)} f\|_{L^2(\mathbb{R})}\lesssim 2^{-\ell} \|f\|_{L^2(\mathbb{R})},
  \]
    where          
 \[
H_{l'}^{(j_1)} f(z) =\int_{\mathbb{R}} f(z - t) e^{i Q(t)} \psi_{l'}^{(j_1)}(t) \Phi_{j_1, j_2}(t) \frac{dt}{t}
  \]
 for $z\in \mathbb{R},$ and $\ell=\ell'-\gamma_k.$ Once the $L^2$ estimate  is proved, one may interpolate with trivial $L^1$ and $L^\infty$ estimates (without any decaying factor $2^{\epsilon \ell}$) to obtain the full range $p\in (1,\infty).$ The decaying factor $2^{\epsilon \ell}$ is also natural to expect: we trivially have 
\[
 \| T_{\gamma_{k_{1}}+l,A}^{\left(k_{1}\right)} f\|_{L^p(\mathbb{R}^n)}\lesssim \|f\|_{L^p(\mathbb{R}^n)}(1<p<\infty),
\]
where $T_{\gamma_{k_{1}}+l,A}^{\left(k_{1}\right)}$ will be defined in \eqref{withA}. However, for the operator $T_{\gamma_{k_{1}}+l,A}^{\left(k_{1}\right)},$ we cannot get the kernel of the operator and use the Fourier transform to estimate the  decay of $2^l$. So,  in this paper, we try to refine the operatot  $T_{\gamma_{k_{1}}+l,A}^{\left(k_{1}\right)}$ with the nature of $A$ to get the following estimate 
\[
 \|T_{\gamma_{k_{1}}+l,A}^{\left(k_{1}\right)}f\|_{L^2(\mathbb{R}^n)}\lesssim 2^{-(\ell-\epsilon)} \|f\|_{L^2(\mathbb{R}^n)}.
\]
The idea is from Chen\cite{Chen2024} and the detail of this part will be discussed in Section 4. 
 
 In what follows, we use the symbol $A\lesssim B$ to denote that there exists a positive constant   $C$ such that $A\le CB.$ 
 
\section{ Preliminaries}

In this section, we recall some facts and lemmas, which we need to prove Theorem \ref{main_th}.

\begin{lemma}\cite{Cohen1981}\label{2.1}
Let $A$ be a function on $\Rn$ with derivatives of order one in $L^q(\mathbb{R}^n)$ for some $q\in (n,\infty]$. Then
\[
|A(x)-A(y)|\leq |x - y|\left(\frac{1}{|I_x^y|}\int_{I_x^y}|\nabla A(z)|^q dz\right)^{1/q},
\]
where $I_x^y$ is the cube which is centered at $x$ and has side length $2|x - y|$.
\end{lemma}

\begin{lemma}\cite{Tao2022}\label{Tao2022}
  Let $\Omega$ be homogeneous of degree zero and $\Omega\in L\log L(S^{n - 1})$, and let $A$ be a function on $\mathbb{R}^n$ with derivatives of order one in $\text{BMO}(\mathbb{R}^n)$. Let $N_{\Omega,A}$ be the maximal operators defined by
\begin{equation}\label{Maxial_Operator1}
    N^*_{\Omega,A}f(x)=\sup_{r > 0}\frac{1}{r^{n + 1}}\int_{|x - y| < r}|\Omega(x - y)||A(x)-A(y)-\nabla A(y)(x - y)||f(y)|dy.
\end{equation}

Then for any $p\in(1,\infty)$,
\[
\|N_{\Omega,A}f\|_{L^p(\mathbb{R}^n)}\lesssim\|\nabla A\|_{\text{BMO}(\mathbb{R}^n)}\|\Omega\|_{L\log L(S^{n - 1})}^*\|f\|_{L^p(\mathbb{R}^n)},
\]
where $\|\Omega\|_{L \log L(S^{n-1})} = \inf\left\{ \lambda > 0 : \int_{S^{n-1}} \frac{|\Omega(x)|}{\lambda} \log\left(e + \frac{|\Omega(x)|}{\lambda}\right) dx \leq 1 \right\}.$
\end{lemma}
\begin{remark}
    It is worth to point out that  $  L(\log L)^2(S^{n - 1})\subset L\log L(S^{n - 1}).$
\end{remark}

\begin{lemma}\cite{169}(Van der Corput lemma ) \label{15}
    Given $\phi$ as a real-valued and smooth function in (a,b), and that $\left|\phi^{(k)}(t)\right| \geq 1$ for all $t \in(a, b)$. Then we have that: 
    \begin{equation*}
        \left|\int_{a}^{b} e^{-i \lambda \phi(t)} \text{d} t\right| \leq C_{k}|\lambda|^{-\frac{1}{k}}, 
    \end{equation*}
    when 
    (i) $k \geq 2$, or
    (ii) $k=1$ and $\phi^{\prime}$ is monotonic,
    the bound $C_k$ is a constant that depends only on k, but not on any a, b, $\phi$, and $\lambda$. 
\end{lemma}

After that, we list some known facts about the Luxemburg norms. Details are given in \cite{Rao1991}. Let \( \Psi : [0, \infty) \to [0, \infty) \) be a Young function, namely, \( \Psi \) is convex and continuous with \( \Psi(0) = 0 \), \( \Psi(t) \to \infty \) as \( t \to \infty \), and \( \Psi(t)/t \to \infty \) as \( t \to \infty \).

We also assume that \( \Psi \) satisfies a doubling condition, i.e., there is \( M \geq 1 \) such that \( \Psi(2t) \leq M\Psi(t) \) for all \( t \in [0, \infty) \).

Let \( \Psi \) be a Young function, and \( Q \subset \mathbb{R}^{n} \) be a cube. Define the Luxemburg norm \( \|\cdot\|_{\mathcal{L}^{\Psi}(Q)} \) by
\[
\|f\|_{\mathcal{L}^{\Psi}(Q)} = \inf\left\{\lambda > 0 : \frac{1}{|Q|}\int_Q \Psi\left(\frac{|f(x)|}{\lambda}\right) dx \leq 1\right\}.
\]
It is well known that
\[
\frac{1}{|Q|}\int_Q \Psi(|f(x)|) dx \leq 1 \Leftrightarrow \|f\|_{\mathcal{L}^{\Psi}(Q)} \leq 1,
\]
and
\begin{align}\label{normequlanr}
    \|f\|_{\mathcal{L}^{\Psi}(Q)} \leq \inf\left\{\mu + \frac{\mu}{|Q|}\int_Q \Psi\left(\frac{|f(x)|}{\mu}\right) dx : \mu > 0\right\} \leq 2\|f\|_{\mathcal{L}^{\Psi}(Q)}. 
\end{align}
 For \( p \in [1, \infty) \) and \( \gamma \in \mathbb{R} \), set \( \Psi_{p, \gamma}(t) = t^p \log^{\gamma}(e + t) \). We denote \( \|f\|_{\mathcal{L}^{\Psi_{p, \gamma}}(Q)} \) as \( \|f\|_{L^{p}(\log L)^{\gamma}, Q} \).

Let \( \Psi \) be a Young function. \( \Psi^* \), the complementary function of \( \Psi \), is defined on \( [0, \infty) \) by
\[
\Psi^*(t) = \sup\{st - \Psi(s) : s \geq 0\}.
\]
The generalization of Hölder inequality
\begin{align}\label{Hlux}
   \frac{1}{|Q|}\int_Q |f(x)h(x)| dx \leq \|f\|_{\mathcal{L}^{\Psi}(Q)} \|h\|_{\mathcal{L}^{\Psi^*}(Q)}
\end{align}
holds for \( f \in \mathcal{L}^{\Psi}(Q) \) and \( h \in \mathcal{L}^{\Psi^*}(Q) \).

For a cube \( Q \subset \mathbb{R}^{n} \) and \( \gamma > 0 \), we also denote \( \|f\|_{\exp L^{\gamma}, Q} \) by
\[
\|f\|_{\exp L^{\gamma}, Q} = \inf\left\{t > 0 : \frac{1}{|Q|}\int_Q \exp\left(\frac{|f(y)|^\gamma}{t}\right) dy \leq 2\right\}.
\]

As it is well known, for \( \Psi(t) = \log(e + t) \), its complementary function \( \Psi^*(t) \approx e^t - 1 \). Let \( b \in \text{BMO}(\mathbb{R}^{n}) \), and the  John-Nirenberg inequality tells us the basic fact  that for any \( Q \subset \mathbb{R}^{n} \) and \( p \in [1, \infty) \),
\[
\|b -\langle b \rangle_Q\|_{\exp L^{1/p}, Q} \lesssim\|b\|_{BMO(\mathbb{R}^{n})}^{p}，
\]
where $\langle b \rangle_Q=\frac{1}{|Q|}\int_{Q}b(x)dx$ and $\|b\|_{BMO}=\sup_{Q}\int_{Q}|b(x)-\langle b \rangle_Q|dx$.

This, together with the inequality \eqref{Hlux}, implies that
\begin{align}\label{Dluxemburg}
    \frac{1}{|Q|} \int_Q |b(x) - \langle b \rangle_Q|^p |h(x)|^p \, dx \lesssim \|h\|_{L^p (\log L), Q}^p \|b\|_{\text{BMO}(\mathbb{R}^{n})}^p.
\end{align}

Now, we turn to prove the main result.
\section{Prove of the Theorem \ref{main_th}}

 \begin{proof}
Equipped with polar coordinates transformation, we can write  as 

\[
T_{Q,A}=\int_{S^{n-1}} \Omega\left(y^{\prime}\right) \int_0^{\infty}\frac{f(x-ty^{\prime})}{t^2}e^{i Q(t)}(A(x)-A(x-ty')-t\nabla A(x-ty')\cdot y') \text{d} t\text{d}\sigma\left(y^{\prime}\right).\] 
To be convince, we denote $A(x,y',t)=A(x)-A(x-ty')-t\nabla A(x-ty')\cdot y'.$
%  Let \(\phi \in C_{0}^{\infty}(\mathbb{R}^{d})\) be a radial function with supp \(\phi \subset B(0,2)\) , and \(\phi(x)=1\) when \(|x| \le1.\) Set \(\varphi(x)=\phi(x)-\phi(2 x).\)  We then have that 
% \[\sum_{j \in \mathbb{Z}} \varphi\left(2^{-j} |x|\right) \equiv 1, \quad|x|>0 .\]

% Let \(\varphi_{j}(|x|)=\varphi(2^{-j} |x|)\) and
% \[
% T_{Q,A}f(x)=\sum_{j\in\mathbb{Z} } T_{Q,A,j}f(x),
% \]
% where
% \[
% T_{Q,A,j}=\int_{S^{n-1}} \Omega\left(y^{\prime}\right) \int_0^{\infty}\frac{f(x-ty^{\prime})}{t^2}e^{i Q(t)}A(x,y',t)\varphi_{j}(t) \text{d} t\text{d}\sigma\left(y^{\prime}\right).\]  
As is shown in the \cite{Guo2017,Ma2024}, we use the fewnomials phases method to estimate this integral. First we split  $R^{+}$ into different intervals and show that there always been a monomial to "dominates" our polynomial $Q.$

Here we assume $1<\alpha_{1}<\cdots<\alpha_{m}$. Denote $\beta$ the degree of polynomial Q, that is $\beta=\alpha_{m}$. Let $\lambda=2^\frac{1}{\beta}$. Define $b_k\in Z$ such that $\lambda^{b_{k}} \leq\left|a_{k}\right|<\lambda^{b_{k}+1}$. We define a few bad scales. For $1 \leq k_{1}<k_{2} \leq m$, define
\begin{equation*}
    \mathcal{K}_{\text {bad }}^{(0)}\left(k_{1}, k_{2}\right):=\left\{l \in \mathbb{Z}: 2^{-C_{0}}\left|a_{k_{2}} \lambda^{\alpha_{k_{2}} l}\right| \leq\left|a_{k_{1}} \lambda^{\alpha_{k_{1}} l}\right| \leq 2^{C_{0}}\left|a_{k_{2}} \lambda^{\alpha_{k_{2}} l}\right|\right\} .
\end{equation*}
Here $C_0$ is a very large number. Notice that $l$ satisfies 
$$-1-\beta C_{0}+b_{k_{2}}-b_{k_{1}} \leq\left(\alpha_{k_{1}}-\alpha_{k_{2}}\right) l \leq \beta C_{0}+b_{k_{2}}-b_{k_{1}}+1.$$
Hence $\mathcal{K}_{\text {bad }}^{(0)}\left(k_{1}, k_{2}\right)$ is a connected set whose cardinality is smaller than $4\beta C_0$. Define 
$$\mathcal{K}_{\text {good }}^{(0)}:=\left(\bigcup_{k_{1} \neq k_{2}} \mathcal{K}_{\text {bad }}^{(0)}\left(k_{1}, k_{2}\right)\right)^{c} .$$
Then, the set $\mathcal{K}_{\text {good }}^{(0)}$ has at most $m^2$ connected components, each of which has a monomial "dominated" Q. Similarly, we define 
\begin{equation*}
    \begin{aligned}
\mathcal{K}_{\text {bad }}^{(1)}\left(k_{1}, k_{2}\right):= & \left\{l \in \mathbb{Z}: 2^{-C_{0}}\left|\alpha_{k_{2}}\left(\alpha_{k_{2}}-1\right) a_{k_{2}} \lambda^{\alpha_{k_{2}} l}\right| \leq\left|\alpha_{k_{1}}\left(\alpha_{k_{1}}-1\right) a_{k_{1}} \lambda^{\alpha_{k_{1}} l}\right|\right. \\
& \left.\leq 2^{C_{0}}\left|\alpha_{k_{2}}\left(\alpha_{k_{2}}-1\right) a_{k_{2}} \lambda^{\alpha_{k_{2}} l}\right|\right\}.
\end{aligned}
\end{equation*}
Moreover,\par
$\mathcal{K}_{\text {bad }}^{(1)}:=\bigcup_{k_{1} \neq k_{2}} \mathcal{K}_{\text {bad }}^{(1)}\left(k_{1}, k_{2}\right)$ and 
$\mathcal{K}_{\text {good }}:=\mathcal{K}_{\text {good }}^{(0)} \backslash \mathcal{K}_{\text {bad }}^{(1)}. $ \par
Analogously, $\mathcal{K}_{\text {good }}$ has at most $m^{4}$ connected components. \par
Firstly, we estimate the bad scales: \par
\textbf{Bad scales: } For every bad scales, there is not a "dominated" monomial to control the polynomial Q. But we can "absorb" the oscillating term using the Maximal operator $N_{\Omega,A}f.$ More precisely speaking, suppose that we are working on the collection of bad scales $\mathcal{K}_{\text {bad }}^{(0)}\left(k_{1}, k_{2}\right)$ for some $k_1$ and $k_2,$ while the same arguments can be given for $\mathcal{K}_{\text {bad }}^{(1)}\left(k_{1}, k_{2}\right).$  Given $\psi_{0}$ be a smooth bump function support on $[\lambda^{-1},\lambda^{2}]$ such that 
$$\sum_{l \in \mathbb{Z}} \psi_{0}\left(\frac{t}{\lambda^{l}}\right)=1, \quad \forall t \neq 0.$$

We denote  \( \psi_{l}(t)= \psi_{0}\left(\frac{t}{\lambda^{l}}\right)\) and

\[
T_{Q,A,l}=\int_{S^{n-1}} \Omega\left(y^{\prime}\right) \int_0^{\infty}\frac{f(x-ty^{\prime})}{t^2}e^{i Q(t)}A(x,y',t) \psi_{l}(t) \text{d} t\text{d}\sigma\left(y^{\prime}\right).\]

Recall the face that the cardinality of $\mathcal{K}_{\text {bad }}^{(0)}\left(k_{1}, k_{2}\right)$ is at most $4\beta C_0$, as \cite{Ma2024}, we could have 
\begin{equation} \label{8}
    \begin{array}{l}
\left|\sum_{l \in \mathcal{K}_{\text {bad }}^{(0)}\left(k_{1}, k_{2}\right)} T_{Q,A,l} f(x)\right| \\
\leq\sum_{j\in \mathcal{K}_{\text {bad }}^{(0)}\left(k_{1}, k_{2}\right)} \left| \int_{S^{n-1}} \Omega\left(y^{\prime}\right)\int_0^{\infty}\frac{f(x-ty^{\prime})}{t^2}e^{i Q(t)}A(x,y',t) \psi_{l}(t)\text{d} t\text{d}\sigma\left(y^{\prime}\right)\right| \\
\leq \sum_{l\in \mathcal{K}_{\text {bad }}^{(0)}\left(k_{1}, k_{2}\right)} \int_{S^{n-1}}\left| \Omega\left(y^{\prime}\right) \right|\int_0^{\infty}\left|\frac{f(x-ty^{\prime})}{t^2}\right|\left|A(x,y',t) \psi_{l}(t)\right|\text{d} t\text{d}\sigma\left(y^{\prime}\right) \\
\leq 4 C_{0} \sum_{l=l_{0}}^{l_{0}+d}\int_{S^{n-1}}\left| \Omega\left(y^{\prime}\right) \right|\int_0^{\infty}\left|\frac{f(x-ty^{\prime})}{t^2}\right|\left|A(x,y',t) \psi_{l}(t)\right|\text{d} t\text{d}\sigma\left(y^{\prime}\right) \\
\le C N_{\Omega, A}^{*} f(x).
\end{array}
\end{equation}
By Lemma \ref{Tao2022}, the estimate of bad scales is completed. 

Now we turn to estimate the good scales,\par
\textbf{Good scales: }For the good scales, suppose we are working on one connected component of $\mathcal{K}_{\text {good }}$. Now for each $l$ in such a component, there always be a "dominate" monomial. We assume that $a_{k_{1}} t^{\alpha_{k_{1}}}$ dominate Q(t) that is, 
$$\left|a_{k_{1}} \lambda^{\alpha_{k_{1}} l}\right| \geq 2^{C_{0}}\left|a_{k_{1}^{\prime}} \lambda^{\alpha_{k_{1}^{\prime}} l}\right|$$ for every $k_{1}^{\prime} \neq k_{1}.$ 

And $a_{k_{2}} \alpha_{k_{2}}\left(\alpha_{k_{2}}-1\right) t^{\alpha_{k_{2}}-2}$ "dominates" $Q^{\prime \prime}(t)$, that is 
$$\left|a_{k_{2}} \alpha_{k_{2}}\left(\alpha_{k_{2}}-1\right) \lambda^{\alpha_{k_{2}} l}\right| \geq 2^{C_{0}}\left|a_{k_{2}^{\prime}} \alpha_{k_{2}^{\prime}}\left(\alpha_{k_{2}^{\prime}}-1\right) \lambda^{\alpha_{k_{2}^{\prime}} l}\right|$$ for every $k_{2}^{\prime} \neq k_{2}.$
We denote such a set $\mathcal{K}_{\text {good }}\left(k_{1}, k_{2}\right). $
Under this assumption, we have the estimates 
$$|Q(t)| \leq 2\left|a_{k_{1}} t^{\alpha_{k_{1}}}\right| \text { and }\left|Q^{\prime \prime}(t)\right| \geq\left|a_{k_{1}} t^{\alpha_{k_{1}}-2}\right|$$
for every $t \in\left[\lambda^{l-2}, \lambda^{l+1}\right]$ with $l \in \mathcal{K}_{\text {good }}\left(k_{1}, k_{2}\right). $We know that $\lambda=2^{1 /\beta}$ is the smallest scale that we will work with. This scale is only visible when $a_{\beta} t^{\beta}$ dominates. When some other monomial dominates, at such a small scale, our polynomial will not have
enough room to see the oscillation. Define $\lambda_{k_{1}}:=2^{1 / \alpha_{k_{1}}},$ we choose this scale because the monomial $a_{k_{1}} t^{\alpha_{k_{1}}}$ dominates. Let 
$$\Phi_{k_{1}, k_{2}}(t)=\sum_{l \in \mathcal{K}_{\text {good }}\left(k_{1}, k_{2}\right)} \psi_{l}(t). $$
Notice that here we join all the small scales from $\mathcal{K}_{\text {good }}\left(k_{1}, k_{2}\right)$ to form a larger scale. Next we will apply a new partition of unity to the function $\Phi_{k_{1}, k_{2}}. $
Let
$$T_{l^{\prime},A}^{\left(k_{1}\right)} f(x)=\int_{S^{n-1}}  \Omega\left(y^{\prime}\right) \int_0^{\infty}\frac{f(x-ty^{\prime})}{t^2}e^{i Q(t)}A(x,y',t) \psi_{l^{\prime}}^{\left(k_{1}\right)}(t) \Phi_{k_{1}, k_{2}}(t)dt \text{d} \sigma\left(y^{\prime}\right),$$
where $\psi_{0}^{\left(k_{1}\right)}$ is a nonnegative smooth bump function supported on $\left[\lambda_{k_{1}}^{-1}, \lambda_{k_{1}}^{2}\right]$ such that
$$\sum_{l^{\prime} \in \mathbb{Z}} \psi_{l^{\prime}}^{\left(k_{1}\right)}(t)=1 \text { for every } t>0 \text {, with } \psi_{l^{\prime}}^{\left(k_{1}\right)}(t):=\psi_{0}^{\left(k_{1}\right)}\left(\frac{t}{\lambda_{k_{1}}^{l^{\prime}}}\right) \text {. }$$
We define $B_{k_{1}} \in \mathbb{Z}$ such that $\lambda_{k_{1}}^{-B_{k_{1}}} \leq\left|a_{k_{1}}\right|<\lambda_{k_{1}}^{-B_{k_{1}}+1},$ and $\gamma_{k_{1}}=B_{k_{1}} / \alpha_{k_{1}}. $  \par
We split the sum in $l^{\prime}$ into two cases.
\begin{equation} \label{14}
    \sum_{l^{\prime} \in \mathbb{Z}} T_{l^{\prime},A}^{\left(k_{1}\right)} f=\sum_{l^{\prime} \leq \gamma_{k_{1}}} T_{l^{\prime},A}^{\left(k_{1}\right)} f+\sum_{l^{\prime}>\gamma_{k_{1}}} T_{l^{\prime},A}^{\left(k_{1}\right)} f.
\end{equation}
For the first sum, $\left|\sum_{l^{\prime} \leq \gamma_{k_{1}}} T_{l^{\prime},A}^{\left(k_{1}\right)} f(x)\right|\le I +II,$
where 
\[I:=\mid \sum_{l^{\prime} \leq \gamma_{k_{1}}}\int_{S^{n-1}}\Omega\left(y^{\prime}\right) \int_0^{\infty}\frac{f(x-ty^{\prime})}{t^2}A(x,y',t)\psi_{l^{\prime}}^{\left(k_{1}\right)}(t) \Phi_{k_{1}, k_{2}}(t)dt \text{d} \sigma\left(y^{\prime}\right)\rvert\, \] and
\[ II:=\sum_{l^{\prime} \leq \gamma_{k_{1}}}\left|\int_{S^{n-1}}\Omega\left(y^{\prime}\right) \int_0^{\infty}\frac{f(x-ty^{\prime})}{t^2}A(x,y',t)\left(e^{i Q(t)}-1\right) \psi_{l^{\prime}}^{\left(k_{1}\right)}(t) \Phi_{k_{1}, k_{2}}(t) \text{d}t  \text{d} \sigma\left(y^{\prime}\right) \right|.\]

First, we consider the term $II.$ We have, by the triangle inequality,
    \begin{align*}
  II\le& \sum_{l^{\prime} \leq \gamma_{k_{1}}}\int_{S^{n-1}}\left|\Omega\left(y^{\prime}\right) \right|\int_0^{\infty}\left|\frac{f(x-ty^{\prime})}{t^2}\right|\\
  &\times\left|A(x,y',t)\right|\left|\left(e^{i Q(t)}-1\right) \right|\psi_{l^{\prime}}^{\left(k_{1}\right)}(t) \Phi_{k_{1}, k_{2}}(t) \text{d}t  \text{d} \sigma\left(y^{\prime}\right) \\
  \le&   \sum_{l^{\prime} \leq \gamma_{k_{1}}}\int_{S^{n-1}}\left|\Omega\left(y^{\prime}\right) \right|\int_0^{\infty}\left|\frac{f(x-ty^{\prime})}{t^2}\right|\\
  &\times\left|A(x,y',t)\right|\left|a_{k_{1}}\right||t|^{\alpha_{k_{1}}-1}\psi_{l^{\prime}}^{\left(k_{1}\right)}(t) \Phi_{k_{1}, k_{2}}(t) \text{d}t  \text{d} \sigma\left(y^{\prime}\right) \\
 \le &  \sum_{l^{\prime} \leq \gamma_{k_{1}}}\int_{S^{n-1}}\left|\Omega\left(y^{\prime}\right) \right|\\
  &\times\int_{\lambda_{k_{1}}^{\gamma_{k_{1}}-l-2}}^{\lambda_{k_{1}}^{\gamma_{k_{1}}-l+1}}\left|\frac{f(x-ty^{\prime})}{t^2}\right|\left|A(x,y',t)\right|\left|a_{k_{1}}\right||t|^{\alpha_{k_{1}}-1}\psi_{l^{\prime}}^{\left(k_{1}\right)}(t) \text{d}t  \text{d} \sigma\left(y^{\prime}\right) \\ 
 \leq&\left|a_{k_{1}}\right|\sum_{l\ge 0} \lambda_{k_{1}}^{\left(\gamma_{k_{1}}-l+1\right)\alpha_{k_{1}}}N_{\Omega, A}^{*} f(x) \le CN_{\Omega, A}^{*} f(x).
      \end{align*}

    Now we consider the $I$. Denote \[z(t)=\sum_{l^{\prime} \leq \gamma_{k_{1}}} \psi_{l^{\prime}}^{\left(k_{1}\right)}(t) \Phi_{k_{1}, k_{2}}(t).\] So $I$ can be controlled by 
    \begin{equation}
        \begin{aligned}
&\left| \int_{{S^{n - 1}}} \Omega  \left( {{y^\prime }} \right)\int_{\left\{ {t \in {\mathbb{R}^ + }:z(t) = 1} \right\}}  \frac{f(x-ty^{\prime})}{t^2}A(x,y',t)\psi_{l^{\prime}}^{\left(k_{1}\right)}(t) \Phi_{k_{1}, k_{2}}(t)dt {\rm{d}}\sigma \left( {{y^\prime }} \right) \right| \\
    &\quad+\int_{S^{n-1}}\left|\Omega\left(y^{\prime}\right)\right| \int_{\left\{t \in \mathbb{R}^{+}: z(t) \neq 1\right\}}\left|\frac{f(x-ty^{\prime})}{t^2}A(x,y',t)\psi_{l^{\prime}}^{\left(k_{1}\right)}(t) \Phi_{k_{1}, k_{2}}(t)\right|dt \text{d} \sigma\left(y^{\prime}\right) \\
    \quad \leq&|\text{p} \cdot \text{v} \cdot \int_{\mathbb{R}^{n}} \frac{\Omega(x-y)}{|x-y|^{d+1}}(A(x)-A(y)-\nabla A(y)(x-y)) f(y)  \text{d} y|\\
    &+\sup _{\varepsilon>0}\left|\int_{|y| \geq \varepsilon}\frac{\Omega(x-y)}{|x-y|^{d+1}}(A(x)-A(y)-\nabla A(y)(x-y)) f(y) \text{d} y\right|\\ 
    &+C\int_{S^{n-1}}|\Omega(y^{\prime})|\int_{\lambda_{k_1}^{\gamma_{k_1}-2}\leq t \leq {\lambda_{k_1}}^{\gamma_{k_1}+1}}\left|\frac{f(x-ty^{\prime})}{t^2}A(x,y',t)\psi_{l^{\prime}}^{\left(k_{1}\right)}(t) \Phi_{k_{1}, k_{2}}(t)dt\right|d\sigma(y^{\prime})\\  
    \quad \leq&\left|T_{\Omega,A} f(x)\right|+T_{\Omega,A}^{*} f(x)+N_{\Omega, A}^{*} f(x).
    \end{aligned}
    \end{equation}
    According to Theorem \ref{TBOUND},  Theorem \ref{MTbound} and Lemma \ref{Tao2022}, we can get the $L^p$ boundedness of $I$.\par
    So we only need to estimate the remained term of \eqref{14}: 
    \[\left\|\sum_{l^{\prime}>\gamma_{k_{1}}} T_{l^{\prime},A}^{\left(k_{1}\right)} f\right\|_{L^{p}(\Rn)} \leq \sum_{l=1}^{\infty}\left\|T_{\gamma_{k_{1}}+l,A}^{\left(k_{1}\right)}f\right\|_{L^{p}(\Rn)}.\]
    It is easy to get that  $ T_{\gamma_{k_{1}}+l,A}^{\left(k_{1}\right)} f$ is controlled by $N_{\Omega, A}^{*} f(x)$. However, we need exponential decay for the $L^p(\forall p\in(1,\infty))$ boundedness of $T_{\gamma_{k_{1}}+l,A}^{\left(k_{1}\right)}.$  Now, we can rewrite  $T_{\gamma_{k_{1}}+l,A}^{\left(k_{1}\right)}$  as   
    \begin{align}\label{withA}
        T_{\gamma_{k_{1}}+l,A}^{\left(k_{1}\right)}f=\int_{\Rn} m^\vee(x-y) A(x,y',t)f(y)dy ,
    \end{align}
where $m(\xi)=\int_{S^{n-1}}\Omega\left(y^{\prime}\right) \int_0^{\infty}\frac{1}{t^2}e^{i (Q(t)+\xi\cdot t y^\prime)}\psi_{l^{\prime}}^{\left(k_{1}\right)}(t) \Phi_{k_{1}, k_{2}}(t) dt \text{d} \sigma\left(y^{\prime}\right).$

Furthermore, it is easy to check that 
       \begin{align*} 
      &\left|\int_{0}^{\infty} e^{i Q(t)+i t y^{\prime} \cdot \xi} \psi_{\gamma_{k_{1}}+l}^{\left(k_{1}\right)}(t) \Phi_{k_{1}, k_{2}}(t) \frac{\text{d}t}{t^2}\right|\\
      &=\lambda_{k_{1}}^{{-(\gamma_{k_{1}}+l)}}\left|\int_{0}^{\infty} e^{i Q\left(\lambda_{k_{1}}^{\gamma_{k_{1}}+l}{ }_{t}\right)+i \lambda_{k_{1}}^{\gamma_{k_{1}}+l}{ }_{t y^{\prime} \cdot \xi}} \psi_{0}^{\left(k_{1}\right)}(t) \Phi_{k_{1}, k_{2}}\left(\lambda_{k_{1}}^{\gamma_{k_{1}}+l} t\right) \frac{\text{d} t}{t^2}\right|. 
    \end{align*}

    We calculate the second order derivative of the phase function: 
    $$\lambda_{k_{1}}^{2 \gamma_{k_{1}}+2 l}\left|Q^{\prime \prime}\left(\lambda_{k_{1}}^{\gamma_{1}+l} t\right)\right| \geq \frac{1}{2}\left|a_{k_{1}}\right| \lambda_{k_{1}}^{B_{k_{1}}+\alpha_{k_{1}} l} \geq 2^{l-2}. $$
    Using the Lemma \ref{15}, we deduce that 
     \begin{align}\label{kerestimat}
         \left|\int_{0}^{\infty} e^{i Q\left(\lambda_{k_{1}}^{k_{1}+l} t\right)+i \lambda_{k_{1}}^{k_{1}+l} t y^{\prime} \cdot \xi} \psi_{0}^{\left(k_{1}\right)}(t) \Phi_{k_{1}, k_{2}}\left(\lambda_{k_{1}}^{\gamma_{k_{1}}+l} t\right) \frac{\text{d} t}{t}\right| \leq C 2^{-\frac{l}{2}}\lambda_{k_{1}}^{-(\gamma_{k_{1}}+l)}.
     \end{align}
So we should find a small enough  $\epsilon\in(0,1)$ such that 
\begin{align}\label{L2}
                     \left\|T_{\gamma_{k_{1}}+l,A}^{\left(k_{1}\right)}(f)\right\|_{L^{2}(\Rn)}\le 2^{-\frac{l}{2}+\epsilon}\left\|f\right\|_{L^{2}(\Rn)},          
            \end{align}
  which will be proved in next section.                                                                            
     Due to the interpolation,  for each $p(1<p<\infty)$, there exists a  $\theta (0<\theta<1)$ such that      
  \begin{align}\label{Lp}
                     \left\|T_{\gamma_{k_{1}}+l,A}^{\left(k_{1}\right)}(f)\right\|_{L^{p}(\Rn)}\le 2^{-\theta(\frac{l}{2}-\epsilon)}\left\|f\right\|_{L^{p}(\Rn)},          
            \end{align}
and  
\[\left\|\sum_{l^{\prime}>\gamma_{k_{1}}} T_{l^{\prime},A}^{\left(k_{1}\right)} f\right\|_{L^{p}(\Rn)} \leq \sum_{l=1}^{\infty}\left\|T_{\gamma_{k_{1}}+l,A}^{\left(k_{1}\right)}f\right\|_{L^{p}(\Rn)}\le \left\|f\right\|_{L^{p}(\Rn)} .\]
Theorem \ref{main_th} is proved.
   
     \end{proof}

\section{the Proof of \eqref{L2}}
   In this section, we also assume $\|\Omega\|_{L\log^2L(S^{n - 1})}^*\le 1.$  We set $E_0=\{x| |\Omega(x)|\le 2\}$,  $E_d=\{x| 2^d\le|\Omega(x)|\le 2^{d+1}\}$ 
and  $\Omega_d=\Omega\chi_{E_d}.$  Let $T_{\gamma_{k_{1}}+l}^{\left(k_{1},d\right)} $   and  $T_{\gamma_{k_{1}}+l,A}^{\left(k_{1},d\right)}$ be the operator defined by
\[
T_{\gamma_{k_{1}}+l}^{\left(k_{1},d\right)}(f)(x)=\int_{\mathbb{R}^d}W^d(x - y)f(y)dy,
\]
and
\[
T_{\gamma_{k_{1}}+l,A}^{\left(k_{1},d\right)}(f)(x)=\int_{\mathbb{R}^d}W^d(x - y)(A(x)-A(y)-\nabla A(y)(x - y))f(y)dy,
\]
where $m^d(\xi)=\int_{S^{n-1}}\Omega_d\left(y^{\prime}\right) \int_0^{\infty}\frac{1}{t^2}e^{i (Q(t)+\xi\cdot t y^\prime)}\psi_{l^{\prime}}^{\left(k_{1}\right)}(t) \Phi_{k_{1}, k_{2}}(t) dt \text{d} \sigma\left(y^{\prime}\right)$ and $W^d=(m^d)^\vee.$

Due to the estimate in \eqref{kerestimat},  it is easy to check that 
\begin{align}\label{3.3}
                     \left\|T_{\gamma_{k_{1}}+l}^{\left(k_{1},d\right)}(f)\right\|_{L^{2}(\Rn)}\lesssim 2^{-\frac{l}{2}}\lambda_{k_{1}}^{-(\gamma_{k_{1}}+l)}\|\Omega_d\|_{L^1(\mathbb{S}^{n-1})}\left\|f\right\|_{L^{2}(\Rn)}.          
            \end{align}
On the other  hand, observe that $\text{supp}\,W^d_j\subset\{x:|x|\leq\lambda_{k_{1}}^{\gamma_{k_{1}}+l}\}$. If $I$ is a cube having side length $\lambda_{k_{1}}^{\gamma_{k_{1}}+l}$, and $f\in L^2(\mathbb{R}^n)$ with $\text{supp}\,f\subset I$, then $T_{\gamma_{k_{1}}+l}^{\left(k_{1},d\right)}f\subset 100nI$. Therefore, to prove (3.14), we may assume that $\text{supp}\,f\subset I$ with $I$ a cube having side length $\lambda_{k_{1}}^{\gamma_{k_{1}}+l}$. Let $x_0$ be a point on the boundary of $200dI$ and $A_I^*(y)=A(y)-\sum_{ j=1}^n \langle\partial_jA\rangle_{ 100nI} y_j$, and
\[
A_I(y)=(A_I^*(y)-A_I^*(x_0))\xi_I(y),
\]
where $\zeta_I\in C_0^{\infty}(\mathbb{R}^n)$, $\text{supp}\,\zeta\subset 150nI$ and $\zeta(x)\equiv1$ when $x\in 100nI$. Observe that $\|\nabla\zeta\|_{L^{\infty}(\mathbb{R}^n)}\lesssim \lambda_{k_{1}}^{\gamma_{k_{1}}+l}$. An application of Lemma \ref{2.1} tells us that for all $y\in 100nI$,
\[
|A_I(y)-A_I(x_0)|\lesssim\lambda_{k_{1}}^{\gamma_{k_{1}}+l}.
\]
This shows that
\[
\|A_I\|_{L^{\infty}(\mathbb{R}^n)}\lesssim\lambda_{k_{1}}^{\gamma_{k_{1}}+l}.
\]
Write
\[
T_{\gamma_{k_{1}}+l,A}^{\left(k_{1}\right)}f(x)=A_I(x)T_{\gamma_{k_{1}}+l}^{\left(k_{1}\right)}f(x)-T_{\gamma_{k_{1}}+l}^{\left(k_{1}\right)}(A_If)(x)-\sum_{j= 1}^{n}[h_j,T_{\gamma_{k_{1}}+l}^{\left(k_{1}\right)}](f\partial_{j}A_I)(x),
\]
where $h_j(x)=x_j.$  It then follows that
\begin{align*}
\|A_IT_{\gamma_{k_{1}}+l}^{\left(k_{1},d\right)}f\|_{L^2(\mathbb{R}^n)}+\|T_{\gamma_{k_{1}}+l}^{\left(k_{1},d\right)}(A_If)\|_{L^2(\mathbb{R}^n)}&\lesssim 2^{-\frac{l}{2}}\|\Omega_d\|_{L^1(\mathbb{S}^{n-1})}\|f\|_{L^2(\mathbb{R}^n)}.
\end{align*}
So we have the following result 
\begin{align*}
&\|A_IT_{\gamma_{k_{1}}+l}^{\left(k_{1}\right)}f\|_{L^2(\mathbb{R}^n)}+\|T_{\gamma_{k_{1}}+l}^{\left(k_{1}\right)}(A_If)\|_{L^2(\mathbb{R}^n)}
\\\le&\sum_{d>0}\left(\|A_IT_{\gamma_{k_{1}}+l}^{\left(k_{1},d\right)}f\|_{L^2(\mathbb{R}^n)}+\|T_{\gamma_{k_{1}}+l}^{\left(k_{1},d\right)}(A_If)\|_{L^2(\mathbb{R}^n)}\right)\\\lesssim& 2^{-\frac{l}{2}}\|\Omega\|_{L^1(\mathbb{S}^{n-1})}\|f\|_{L^2(\mathbb{R}^n)}.
\end{align*}

Applying the John - Nirenberg inequality, we know that
\[
\|\partial_{j}A_I\|_{\text{exp}L^{1/2},I}^2\lesssim 1.
\]
Recall that $\text{supp}\,[h_j,T_{\gamma_{k_{1}}+l}^{\left(k_{1},d\right)}](f\partial_{j}A_I)\subset 100nI$. It follows from inequality (2.4) that
\begin{align*}
\|[h_j,T_{\gamma_{k_{1}}+l}^{\left(k_{1}\right)}](f\partial_{j}A_I)\|_{L^2(\mathbb{R}^n)}&=\sup_{\|g\|_{L^2(\mathbb{R}^n)}\leq1}\left|\int_{\mathbb{R}^d}\partial_{j}A_I(x)f(x)[h_j,T_{\gamma_{k_{1}}+l}^{\left(k_{1}\right)}]g(x)dx\right|\\
&\leq\|f\|_{L^2(\mathbb{R}^n)}\sup_{\substack{\|g\|_{L^2(\mathbb{R}^n)}\leq1\\ \text{supp}\,g\subset 100nI}}\|\partial_{j}A_I[h_j,T_{\gamma_{k_{1}}+l}^{\left(k_{1}\right)}]g\|_{L^2(I)}\\
&\leq\|f\|_{L^2(\mathbb{R}^n)}\left(|I|\sup_{\substack{\|g\|_{L^2(\mathbb{R}^n)}\leq1\\ \text{supp}\,g\subset 100nI}}\|[h,T_{\gamma_{k_{1}}+l}^{\left(k_{1}\right)}]g\|_{L^2(\log L)^2,I}^2\right)^{1/2}.
\end{align*}
Now let $g\in L^2(\mathbb{R}^n)$ with $\|g\|_{L^2(\mathbb{R}^n)}\leq1$ and $\text{supp}\,g\subset 100nI$. Observe that

\[
\|T_{\gamma_{k_{1}}+l}^{\left(k_{1},d\right)}g\|_{L^{\infty}(\mathbb{R}^n)}\lesssim\|\Omega_d\|_{L^\infty(\mathbb{S}^{n-1})}\lambda_{k_{1}}^{-n(\gamma_{k_{1}}+l)}\|g\|_{L^1(\mathbb{R}^n)},
\]

and so
\[
\|\chi_I[h_j,T_{\gamma_{k_{1}}+l}^{\left(k_{1},d\right)}]g\|_{L^{\infty}(\mathbb{R}^n)}\lesssim \|g\|_{L^1(\mathbb{R}^n)}\lesssim \|\Omega_d\|_{L^\infty(\mathbb{S}^{n-1})}\lambda_{k_{1}}^{-\frac{n(\gamma_{k_{1}}+l)}{2}},
\]
since $supp(g)\subset 100nI$ and 
\[
\|g\|_{L^1(\mathbb{R}^n)}\lesssim|I|^{1/2}\|g\|_{L^2(\mathbb{R}^n)}\lesssim \lambda_{k_{1}}^{\frac{n(\gamma_{k_{1}}+l)}{2}}.
\]
On the other hand, we deduce from \eqref{3.3} that
\begin{align*}
\|[h_j,T_{\gamma_{k_{1}}+l}^{\left(k_{1},d\right)}]g\|_{L^2(\mathbb{R}^n)}&\lesssim \lambda_{k_{1}}^{\gamma_{k_{1}}+l}\|T_{\gamma_{k_{1}}+l}^{\left(k_{1},d\right)}g\|_{L^2(\mathbb{R}^n)}\\
&\lesssim  2^{-\frac{l}{2}}\|\Omega_d\|_{L^1(\mathbb{S}^{n-1})}\|g\|_{L^2(\mathbb{R}^n)}.
\end{align*}
Set $\lambda_0 =  2^{-l}\lambda_{k_{1}}^{-n(\gamma_{k_{1}}+l)} \log^2\left(e+2^{\frac{l}{2}}\right) .$ A straightforward computation tells us that
\begin{align*}
&\int_{I}|[h_j,T_{\gamma_{k_{1}}+l}^{\left(k_{1}\right)}]g(x)|^2\log^2\left(e+\frac{|[h_j,T_{\gamma_{k_{1}}+l}^{\left(k_{1}\right)}]g(x)|}{\sqrt{\lambda_0}}\right)dx\\
\le& \sum_{d>0}\int_{I}|[h_j,T_{\gamma_{k_{1}}+l}^{\left(k_{1},d\right)}]g(x)|^2\log^2\left(e+\frac{|[h_j,T_{\gamma_{k_{1}}+l}^{\left(k_{1},d\right)}]g(x)|}{\sqrt{\lambda_0}}\right)dx\\
\lesssim&\sum_{d\ge0}\log^2\left(e+\frac{\|[h_j,T_{\gamma_{k_{1}}+l}^{\left(k_{1},d\right)}]g(x)\|_{L^{\infty}(\mathbb{R}^n)}}{\sqrt{\lambda_0}}\right)\int_{I}|[h_j,T_{\gamma_{k_{1}}+l}^{\left(k_{1},d\right)}]g(x)|^2dx\\
\lesssim &\sum_{d\ge0}\log^2\left(e+\|\Omega_d\|_{L^\infty(\mathbb{S}^{n-1})}2^\frac{l}{2}\right) 2^{-\frac{l}{2}}\|\Omega_d\|^2_{L^1(\mathbb{S}^{n-1})}\\
\lesssim &\lambda_0\lambda_{k_{1}}^{n(\gamma_{k_{1}}+l)}
\end{align*}
due to the fact that$\frac{|[h_j,T_{\gamma_{k_{1}}+l}^{\left(k_{1}\right)}]g(x)|}{\sqrt{\lambda_0}}\le  2^\frac{l}{2} 
$ and \cite[Lemma 3]{Hu2003}. So we get that there exist an  small enough  $\epsilon$ such that     
\[\|[h_j,T_{\gamma_{k_{1}}+l}^{\left(k_{1}\right)}](f\partial_{j}A_I)\|_{L^2(\mathbb{R}^n)}\lesssim  2^{-l}\log^2\left(e+2^{\frac{l}{2}}\right)  \|f\|_{L^2(\mathbb{R}^n)}\lesssim  2^{-l+\epsilon}  \|f\|_{L^2(\mathbb{R}^n)},\]
 and   inequality \eqref{L2} is proved.


\end{document}